\documentclass[12pt,reqno]{amsart}

\usepackage{verbatim}
\usepackage{amssymb}
\usepackage[dvipsnames]{xcolor}
\usepackage[pdftex]{hyperref}

\hypersetup{colorlinks=true, linkcolor=blue, citecolor=ForestGreen}

\usepackage{amsmath}
\usepackage{amssymb}

\numberwithin{equation}{section}

\author{Ernie Croot}
\address[Ernie Croot]%
  {Department of Mathematics, Georgia Institute of Technology\\
Atlanta, Georgia, U.S.A., 30332.}
\email{ecroot@math.gatech.edu}

\author{Vsevolod F. Lev}
\address[Vsevolod F. Lev]{Department of Mathematics, University of Haifa 
  at Oranim, Tivon, Israel.}
\email{seva@math.haifa.ac.il}

\author{P\'eter P\'al Pach}
\address[P\'eter P\'al Pach]{Department of Computer Science and Information 
  Theory,  Budapest University of Technology and Economics\\ 
M\H{u}\-e\-gye\-tem rkp. 3., H-1111 Budapest, Hungary;\\
MTA-BME Lend\"ulet Arithmetic Combinatorics Research Group\\
  M\H{u}egyetem rkp. 3., H-1111 Budapest, Hungary.}
\email{pach.peter@vik.bme.hu}

\title[Past and future of the cap set problem]%
  {Past and future of the cap set problem}

\begin{document}
\baselineskip=16pt

\maketitle

\tableofcontents

\section{Roth's problem on three-term progressions}

In 1953 K. F. Roth \cite{roth} proved that the largest subset of $[N] := 
\{1,2,...,N\}$ containing no three-term arithmetic progression $x, x+d, x+2d$ 
has size $o(N)$.  Working through his proof (suitably interpreted) one can 
even get a quantitative bound of the form $O(N / \log\log N)$.  This then 
naturally leads to the following question. 
\bigskip

\noindent {\bf Roth's Problem.} What is the size of the largest subset $S 
\subseteq [N]$ containing no three-term arithmetic progressions? 
\bigskip

Most of the progress on this problem since Roth's seminal work makes heavy 
use of a ``density increment argument" pioneered by him.  The idea is that if 
one assumes $S \subseteq [N]$ has no three-term progression, and if $|S| = 
\alpha N$ and $N > N_0(\alpha)$, then one can show that there exists an 
arithmetic progression $P := \{a, a+d, a+2d, ..., a+kd\} \subseteq [N]$, 
where $k > N^{1/2-o(1)}$, such that $|S \cap P| \geq \alpha(1+ c \alpha)|P|$ 
for some $c > 0$.  By translating and rescaling, one then has a 
progression-free set $S' \subseteq [N']$, $|N'| > N^{1/2-o(1)}$, $|S'| \geq 
\alpha(1 + c \alpha) |N'|$. Iterating this (staying above the $N_0(\alpha)$ 
threshold for the interval length), eventually one reaches a contradiction if 
$\alpha > c' /\log\log N$, because if $\alpha$ is this big, one of the sets 
$S''$ so constructed would have to have density $1$, yet also is 
progression-free.  Thus, if the original $S$ is progression-free, then $|S| 
\ll N / \log\log N$. 

Further refinements on the idea included achieving a greater density 
increment per iteration relative to the length of the interval 
\cite{heathbrown, szemeredi}, resulting in bounds for progression-free sets 
of the type $|S| < N (\log N)^{-\delta}$, for some $0 < \delta < 1/2$.  
Replacing density-increments on sub-progressions (as in Roth's method) with 
density-increments on so-called Bohr-neighborhoods, Bourgain \cite{bourgain} 
achieved a bound of the form $|S| \ll N \sqrt{\log \log N \over \log N}$. 
Then in a series of papers by himself \cite{bourgain2} and Sanders 
\cite{sanders, sanders2} the bound was improved to $|S| < N (\log 
N)^{-1+o(1)}$.  Improving this bound even a little bit (lowering the $-1$ to 
$-1-\varepsilon$) would establish the special case $k=3$ of the following 
famous conjecture \cite{erdos}, which if proved would give a far-reaching 
generalization of Szemer\'edi's Theorem \cite{szemeredi2}. 
\bigskip

\noindent {\bf Erd\H os-Tur\'an Conjecture on $k$-term arithmetic 
progressions.}  If $A$ is a set of positive integers such that $\sum_{a \in 
A} 1/a$ diverges, then for every $k \geq 2$, the set $A$ contains a $k$-term 
arithmetic progression. 
\bigskip

The best quantitative bounds in the direction of addressing this theorem in 
the general case (for all values of $k$) are due to Leng, Sah, and Sawhney 
\cite{leng}, who recently proved that for every $k \geq 5$ there exists $c_k 
> 0$ such that the largest subset $S$ of $[N]$ having no $k$-term arithmetic 
progressions has size $|S|\ \ll\ N \exp( - (\log\log N)^{c_k})$. This 
improved upon Gowers's bounds \cite{gowers} that $|S| \ll N (\log \log 
N)^{-2^{-2^{k+9}}}$.  In the case $k=4$, Green and Tao \cite{greentao1, 
greentao2} established the bound $|S| \ll N (\log N)^{-c}$ for some $0 < c < 
1$. 

Bloom and Sisask \cite{bloomsisask} were the first to prove the above 
conjecture for $k=3$, building on the work of Bateman and Katz 
\cite{batemankatz}, by showing that for $N > N_0$ the largest 
progression-free set $S \subseteq [N]$ has size $|S| < N (\log 
N)^{-1-\varepsilon}$ (for some explicit $\varepsilon > 0$).  Then, in a 
remarkable breakthrough, Kelley and Meka \cite{kelleymeka, bloomsisask2} 
improved this to $|S| < N \exp(-c (\log N)^{1/11})$, which Bloom and Sisask 
\cite{bloomsisask3} refined to give $|S| < N \exp(-c' (\log N)^{1/9})$.  
These bounds are not far off from the best possible, since from the work of 
Behrend \cite{behrend} it was known that there exists a three-term 
progression-free set $S \subseteq [N]$ satisfying $|S| > N \exp(-(2\sqrt{\log 
4}+o(1))\sqrt{\log N})$.  This was improved by Elkin \cite{elkin} by a small 
factor tending to infinity, and then recently Elsholtz, Hunter, Proske, and 
Sauermann \cite{elsholtzhunter} gave a substantial further improvement $|S| > 
N \exp(-(C+o(1)) \sqrt{\log N})$, where $C = 2 \sqrt{\log(24/7) \log(2)} < 2 
\sqrt{\log 4}$. 

\section{Finite field settings}

As we saw, the main difficulty in Roth's original approach was getting a high 
enough density increment of the set along progressions, relative to their 
(the progressions) size.  Meshulam \cite{meshulam} considered what this 
argument gives in the case where instead of working with subsets of intervals 
in the integers, one works with subsets of the finite field vector space 
${\mathbb F}_p^n$. The case $p=3$ is known as the {\it cap set problem}.  

In Meshulam's treatment of the general case ${\mathbb F}_p^n$, rather than 
getting a density increment inside a  sub-progression at each iteration (of 
Roth's argument), one gets a density increment on affine subspaces 
(translates of subspaces) $t + V$ where ${\rm dim}(V) = n-1$.  Since these 
affine subspaces are $p^{n-1}$ in size, one can run the density increment 
argument for more steps than if one's sets $S$ were drawn from integer 
intervals $[N]$ when $N \approx p^n$; and, furthermore, the whole argument is 
more elegant and simpler than the integer case, while also containing many of 
the same, or analogous, difficulties. In fact, this is  true of many additive 
combinatorial problems \cite{green, wolf, peluse}. Thus, it is often fruitful 
when trying to solve a problem over ${\mathbb Z}$, say, to first see what one 
can prove for an ${\mathbb F}_p^n$ analogue of that problem. 

In the end, Meshulam proved that the largest subset $S \subseteq {\mathbb 
F}_p^n$ without three-term progressions (or solutions to $x+y = 2z$) 
satisfies 
$$
|S|\ <\ {c_p p^n \over n}.
$$
Meshulam's proof uses Fourier methods, but in \cite{lev2} Lev developed a 
purely combinatorial approach to achieve the same bounds. 

Significantly improving upon Meshulam's bound was considered a major 
challenge, and Terry Tao \cite{tao} even once referred to the overall problem 
of understanding the size of sets without three-term progressions in 
${\mathbb F}_3^n$ as ``perhaps my favorite open question".  

Bateman and Katz were the first to make major progress on it, proving that 
there exists $\varepsilon > 0$ such that in the case $p=3$ one has $|S| \ll 
3^n / n^{1+\varepsilon}$.  Then Ellenberg and Gijswijt \cite{ellenbergG}, 
building on our work in \cite{crootlevpach}, used algebraic methods to prove 
that for every prime $p \geq 3$ there exists $\delta_p > 0$ such that $|S| 
\ll_p (p - \delta_p)^n$.  Further algebraic generalizations of the method 
were given by Tao, Sawin \cite{tao2, tao3}, and Petrov \cite{petrov}.  

More recently, Kelley and Meka \cite{kelleymeka} have developed a 
combinatorial argument (one ingredient of which being \cite{crootsisask}) to 
prove weaker bounds, but still much stronger than other combinatorial and 
Fourier-analytic approaches, achieving $|S| \ll 2^{-\kappa_p n^{1/9}}p^n$.   

Lower bounds were proved by Edel \cite{edel} for $p=3$ giving the existence 
of a set $S$ without three-term progressions that satisfies $|S| > 
(2.217389)^n$.  This then was improved upon by Tyrrell \cite{tyrrell} to $|S| 
> (2.218)^n$, by Romera-Paredes {\it et al} \cite{paredes} to $|S| > (2.2202)^n$, and by 
Naslund \cite{naslund2} to $|S| > (2.2208)^n$.  Recently, Elsholtz, Hunter, Proske, and Sauermann 
\cite{elsholtzhunter} achieved a general lower bound of the shape $|S| > (c 
p)^n$ for some $c > 1/2$ for all primes $p \geq 3$. 

\section{The rise of algebraic methods}

Algebraic methods have been used in several ways in the fields of finite 
geometries, additive combinatorics, and additive number theory.  For example, 
the Chevalley-Warning theorem can be used to quickly prove a special case of 
Olson's theorem \cite{alon} (among many other uses of it); and Stepanov's 
method \cite{stepanov} can be used to count points on curves over a finite 
field.    

More relevant to our discussion is perhaps Alon's Combinatorial 
Nullstellensatz \cite{alon2}, one version of which is: 
\bigskip

\noindent {\bf Theorem 1.} Suppose ${\mathbb F}$ is a field and let 
$f(x_1,...,x_n) \in F[x_1,...,x_n]$.  Suppose the degree deg$(f)$ of $f$ is 
$\sum_{i=1}^n t_i$, where each $t_i$ is a non-negative integer, and suppose 
the coefficient of $\prod_{i=1}^n x_i^{t_i}$ in $f$ is nonzero.  Then if 
$S_1,...,S_n$ are subsets of $F$ with $|S_i| \geq t_i+1$, there are $s_1 \in 
S_1$, $s_2 \in S_2,\dotsc,s_n \in S_n$ so that $f(s_1,...,s_n) \neq 0$. 
 \bigskip
 
To apply this sort of result, one needs to encode the combinatorial problem 
under consideration in terms of vanishing of some low-degree polynomial, and 
then show that the properties of the polynomial (reflecting the original 
combinatorial problem) are inconsistent with the low-degree condition.

This style of reasoning was used in our paper \cite{crootlevpach} on 
three-term progressions in ${\mathbb Z}_4^n$, as we will now discuss. 

In \cite{lev} Lev had generalized Meshulam's result to arbitrary finite 
additive abelian groups $G$, showing that if $S \subseteq G$ has no 
three-term progressions then $|S| < 2 |G| / {\rm rank}(2G)$. Here,  $2G = 
\{2g : g \in G\}$, and ${\rm rank}(H)$ denotes the unique number $r$ in a 
decomposition $H \cong {\mathbb Z}_{d_1} \oplus \cdots \oplus {\mathbb 
Z}_{d_r}$, $d_1 | d_2 | \cdots | d_r$.  Note that in the case $G = {\mathbb 
Z}_4^n$ we have that ${\rm rank}(2G) = n$, so $|S| < 2\cdot 4^n/n$.  Using 
Fourier methods, Tom Sanders \cite{sanders3} gave a stronger bound $|S| = 
o(4^n/n)$. And in our work \cite{crootlevpach} we used algebraic methods to 
prove 
$$
|S|\ <\ 4^{cn},\ {\rm where\ }c \approx 0.926,
$$
thus giving an ``exponential improvement" over previous results. 

A key lemma in our work was the following. 
\bigskip

\noindent {\bf Lemma 1.}  Let ${\mathbb F}$ be a field. Suppose $n \geq 1$ 
and $d \geq 0$ are integers, and let $f \in {\mathbb F}[x_1,...,x_n]$ be a 
multilinear polynomial (that is, all monomials are square-free) of degree at 
most $d$.  Suppose $A \subseteq {\mathbb F}^n$ satisfies $|A| > 2 \sum_{0 
\leq i \leq d/2} {n \choose i}$.  If $f(a-b) = 0$ for every $a,b \in A$ with 
$a \neq b$, then $f(0) = 0$. 
\bigskip

The way this lemma can be used to deduce strong bounds on progression-free 
sets in ${\mathbb Z}_4^n$ is as follows. (We will not give here an optimized 
version of the argument with the bounds claimed above, but rather just an 
easy-to-follow one.) 

First, let $F_n \leq {\mathbb Z}_4^n$ be the subgroup of the $2^n$ elements 
of $\{0,2\}^n$, and for a subset $A \subseteq {\mathbb Z}_4^n$ we let $A_t$ 
denote the set $F_n\cap (A-t)$.  We note that $A_t$, $A_t + A_t$ and $2*A = 
\{2a\ :\ a \in A\}$ all are subsets of $F_n$.  As an additive group we have 
that $F_n$ is isomorphic to ${\mathbb F}_2^n$; and so we can treat these 
three sets as subsets of ${\mathbb F}_2^n$.  In fact, if $\rho : F_n \to 
{\mathbb F}_2^n$ is the obvious group isormophism, we can define $A' = 
\rho(2*A)$ and define $A'_t = \rho(A_t)$.  

Now we suppose we have a set $S\subseteq {\mathbb Z}_4^n$ having no 
three-term progressions.  And let us suppose, for simplicity of discussion, 
that each of the sets $S_t$ either has $0$ elements or has $N$ elements, for 
some $N$ (one could imagine removing some elements from $S$ until this 
``either empty or $N$ elements" condition holds, without much shrinking the 
size of $S$).  So, $|S'| = |\rho(2*S)| = |S|/N$.  

The set $S$ having no three-term progressions implies that all the restricted 
sumsets $S_t \hat + S_t + 2t = 
 \{s_1 + s_2 + 2t\ :\ s_1, s_2 \in S_t,\ s_1 \neq s_2\}$ are disjoint from $2*S$.  The same will be true of
 $S'_t \hat + S'_t + \rho(2t) \subseteq {\mathbb F}_2^n$ and 
 $S' = \rho(2*S) \subseteq {\mathbb F}_2^n$. 

The idea now is to let $f(x_1,...,x_n) \in {\mathbb F}_2[x_1,...,x_n]$ be a 
multilinear polynomial of as low a degree as possible that vanishes on 
$\overline{S'} = {\mathbb F}_2^n \setminus S'$.  Given a degree $d$ we know 
there are $\sum_{i=0}^d {n \choose i}$ square-free monomials in $x_1,...,x_n$ 
of degree at most $d$; and an easy degrees-of-freedom or dimension-counting 
argument shows that if this sum exceeds $|\overline{S'}| = 2^n - |S|/N$, then 
there exists such a polynomial (that vanishes on $\overline{S'}$) of degree 
at most $d$.  Furthermore, this polynomial $f$ will be non-zero and does not 
vanish on all ${\mathbb F}_2^n$. 

We will assume $d$ is minimal such that this holds.  

Now, for every $t$ such that $S'_t \neq \emptyset$, since $S'_t \hat + S'_t + 
\rho(2t) \subseteq \overline{S'}$, we would have the polynomial 
$g(x_1,...,x_n) = f((x_1,...,x_n) + \rho(2t))$ vanishes on $S'_t \hat + 
S'_t$.  By Lemma 1 if $|S'_t| = N > 2 \sum_{0 \leq i\leq d/2} {n \choose i}$, 
then we would also have that $g((0,...,0)) = 0$, which means $f(\rho(2t)) = 
0$.  Since this holds for all those $t$ with $S'_t \neq \emptyset$ it would 
follow that $f$ also vanishes on $S'$.  Since $f$ vanishes on $\overline{S'}$ 
and $S'$, $f$ vanishes on all of ${\mathbb F}_2^n$, which is a contradiction.  

If one now considers the possibilities for $|S|$ and $N > |S|/2^n$ so that 
both $N \leq 2 \sum_{0 \leq i \leq d/2} {n \choose i}$ and $\sum_{i =0}^d {n 
\choose i} > 2^n - |S|/N$ hold, one will see this forces $|S| < 4^{cn}$ for 
some $c > 0$. 
\bigskip

Ellenberg and Gijswijt \cite{ellenbergG} adapted the algebraic argument in 
the ${\mathbb Z}_4^n$ case to prove similar bounds for ${\mathbb F}_p^n$.  
Their proof turned out to be simpler, partly because they didn't have to deal 
with an analogue of cosets of $F_n$.   

They proved that for every prime $p \geq 3$ there exists $0 < c_p < 1$ such 
that if $S \subseteq {\mathbb F}_p^n$ contains no three-term progressions, 
then 
$$
|S|\ \ll\ p^{c_p n}.
$$
\bigskip

Taking inspiration from these papers, Terry Tao \cite{tao2, tao3} introduced 
what he called a ``symmetric formulation" of the methods from 
\cite{crootlevpach} and \cite{ellenbergG}.  He and Will Sawin \cite{tao3} 
introduced the so-called {\it slice-rank}, which for the case of $3$ 
variables $x,y,z$ (the case of interest to proving bounds on sets without 
three-term progressions) can be defined as follows. 
\bigskip

\noindent {\bf Slice-Rank.}  Suppose that ${\mathbb F}$ is a field, 
$A\subseteq\mathbb F$ is a finite set, and $f$ is an $\mathbb F$-valued 
function on the cross product $A\times A\times A$. The slice-rank of $f$ is 
the minimum number $d \geq 1$ such that one can write $f$ as a linear 
combination (over $\mathbb F$) of functions of the forms $g_1(x)h_1(y,z)$, 
$g_2(y) h_2(x,z)$, and $g_3(z)h_3(x,y)$. 
\bigskip

And one of the results he proved about this is the following. 
\bigskip

\noindent {\bf Lemma 2.}  Suppose $A$ is a finite subset of a field ${\mathbb 
F}$ and suppose that $f(x,y,z) : A \times A \times A \to {\mathbb F}$ is the 
``diagonal map" -- that is, $f(x,y,z) = 1$ if $x=y=z$, and is $0$ otherwise.  
Then the slice-rank of $f$ is $|A|$. 
\bigskip

The idea for how to apply this is to assume $S \subseteq {\mathbb F}_3^n$, 
say, has no three-term progressions.  Then, $f(\vec x, \vec y,\vec z) = 
\prod_{i=1}^n (1 - (x_i + y_i + z_i)^2)$ is the diagonal map on $S \times S 
\times S$, because first note that $f$ takes either the value $0$ or $1$ (it 
cannot take the value $-1$); and then in order to be $1$ we would have to 
have all $x_i + y_i + z_i = 0$, which would mean $x+y+z=0$.  Then, since $S$ 
has no three-term progressions, this could only happen if $x=y=z$.  

Next, we expand $f$ into monomials $x_1^{i_1} \cdots x_n^{i_n} 
y_1^{j_1}\cdots y_n^{j_n} z_1^{k_1}\cdots z_n^{k_n}$, where the exponents are 
in $\{0,1,2\}$ and have sum $\leq 2n$ (since the degree of $f$ is $2n$).  And 
now the idea is to write this linear combination of monomials as a linear 
combination of functions of the form $f(x)g(y,z)$, $f(y)g(x,z)$, and of the 
form $f(z)g(x,y)$.  For each choice of $i_1,...,i_n$ with $i_1 + \cdots + i_n \leq 2n/3$ we group all the 
$y_\ell$'s and $z_m$'s together that appear and call that $g(y,z)$, and then $f(x) = 
x_1^{i_1}\cdots x_n^{i_n}$.  We do a similar grouping for each choice of $j_1,...,j_n$ when $j_1+\cdots +j_n \leq 2n/3$ for all the remaining terms (after excluding those where $i_1 +\cdots + i_n \leq 2n/3$ we already counted) in the monomial expansion of $f$, except we get functions of the form $f(y)g(x,z)$; and then the remaining terms will all have $k_1+\cdots + k_n \leq 2n/3$, and we get functions of the form $f(z)g(x,y)$. 

If one counts up the number of different functions of each of the three types 
($f(x)g(y,z)$ or $f(y)g(x,z)$ or $f(z)g(x,y)$), one gets a linear combination 
involving at most 
$$
3 \sum_{a+b+c = n,\ b+2c \leq 2n/3 \atop a,b,c \geq 0}
{n! \over a! b! c!}.
$$
terms.  Since this is an upper-bound on the slice-rank of $f$, Lemma 2 above 
tells us it is also an upper bound on $|S|$.  And now it is not hard to see 
that this upper bound has the form $3^{\kappa n}$ for some $0 < \kappa < 1$. 

\section{Further applications}

One of the early applications of the various methods \cite{crootlevpach, 
ellenbergG, tao3} from the previous section was the work of Naslund and Sawin 
\cite{naslundsawin} on upper bounds for $3$-sunflower-free sets.  That is, 
suppose ${\mathcal F}$ is a family of subsets of $\{1,2,...,n\}$ that does 
not contain a triple of sets $A,B,C$ with the property $A\cap B = A \cap C = 
B \cap C$.  

Erd\H os and Szemer\'edi \cite{erdossz} proved that that any such family 
${\mathcal F}$ must satisfy $|{\mathcal F}| < 2^n \exp(-c\sqrt{n})$.  Then, 
Alon, Shpilka, and Umans \cite{alon3} showed that upper bounds on the size of 
capsets (progression-free sets in ${\mathbb F}_3^n$) translate into upper 
bounds on the size of $3$-sunflower-free sets; and then using the capset 
bounds from \cite{ellenbergG} one obtains a bound of the shape $|{\mathcal 
F}| < c^n$, for some $c = \sqrt{1 + 2.7552} = 1.9378...$.  However, Naslund 
and Sawin \cite{naslundsawin} further strengthened this by applying the 
polynomial method {\it directly} to the problem (rather than passing through 
capset bounds) to obtain the stronger bound $|{\mathcal F}| < 
(2/2^{2/3})^{n(1+o(1))} \approx 1.889881574^{n(1+o(1))}$.   

In \cite{blasiak} Blasiak, Cohn, Grochow, Naslund, Sawin, and Umans used 
these algebraic methods to rule out the existence of a certain type of fast 
matrix-multiplication algorithm that could multiply two $n \times n$ matrices 
in time $n^{2+o(1)}$.  This type of algorithm had been conjectured to exist 
by Cohn, Kleinberg, Szegedy, and Umans \cite{cohn}. 

\section{Directions}

Here we list a few questions worthy of consideration. 
\begin{itemize}
\item Can algebraic methods be used to estimate the size of the largest set 
    $S$ without a $k$-term progression in ${\mathbb F}_p^n$, for $k \geq 
    4$?  

\item Can one use the algebraic methods in the restricted difference 
    settings? For example, how large can a subset $S\subset{\mathbb F}_3^n$ 
    be given that $S$ does not contain any three-term arithmetic 
    progression with the difference in $\{0,1\}^n$? 

\item Along the same lines -- but more general -- is the question of 
    addressing exactly which problems can be solved by the kinds of 
    algebraic methods in this paper.  The work \cite{petrov} is perhaps a 
    path towards addressing this. 

\item Can these methods be extended somehow to address questions in the 
    integers?  Progress on three-term progressions in subsets of integer 
    intervals is already fairly advanced, thanks to the recent work of 
    Kelley and Meka \cite{kelleymeka}; however, it would be nice to have 
    other approaches. 

\item Is there a way to unify all the different algebraic methods for 
    proving combinatorial statements, such as uses of Chevalley-Warning, 
    Stepanov's method, Alon's Nullstellensatz, and now the methods for 
    proving strong bounds on sets without three-term arithmetic 
    progressions?  Are they all really ``the same method" in some sense? 
\end{itemize}

\vfill

\end{document}